\documentclass[a4paper]{amsart}
\usepackage{amssymb}
\usepackage[]{graphicx}
\title{Extremal spectral properties of Otsuki tori}
\author{Alexei V. Penskoi}
\address{Department of Geometry and Topology, 
Faculty of Mathematics and Mechanics, Moscow State University,
Leninskie Gory, GSP-1, 119991, Moscow, Russia\newline \emph{and}
\newline Independent University of Moscow, 
Bolshoy Vlasyevskiy pereulok~11, 119002, Moscow, Russia}
\email{penskoi@mccme.ru}
\date{}
\newtheorem{Proposition}{Proposition}
\newtheorem{Definition}{Definition}
\newtheorem*{Theorem}{Theorem}

\DeclareMathOperator{\Area}{Area}
\DeclareMathOperator{\SO}{SO}
\begin{document}
\begin{abstract}
Otsuki tori form a countable family
of immersed minimal two-dimensional
tori in the unitary three-dimensional sphere.
According to El~Soufi-Ilias theorem,
the metrics on the Otsuki tori are extremal
for some unknown eigenvalues of the 
Laplace-Beltrami operator. Despite
the fact that the Otsuki tori are defined
in quite an implicit way, we
find explicitly the numbers of 
the corresponding extremal eigenvalues.
In particular we provide an extremal
metric for the third eigenvalue
of the torus.
\end{abstract}
\maketitle

\section*{Introduction}

Let $M$ be a closed surface and $g$ be
a Riemannian metric on $M.$
Let us consider the associated Laplace-Beltrami
operator $\Delta:C^\infty(M)\longrightarrow C^\infty(M),$
$$
\Delta f=-\frac{1}{\sqrt{|g|}}\frac{\partial}{\partial x^i}%
\left(\sqrt{|g|}g^{ij}\frac{\partial f}{\partial x^j}\right).
$$
It is well-known that the eigenvalues
\begin{equation}\label{eigenvalues}
0=\lambda_0(M,g)<\lambda_1(M,g)\leqslant%
\lambda_2(M,g)\leqslant\lambda_3(M,g)\leqslant\dots
\end{equation}
of $\Delta$ possess the following rescaling property,
$$
\forall t>0\quad\lambda_i(M,tg)=\frac{\lambda_i(M,g)}{t}.
$$
Hence, given a fixed surface $M,$ we have 
$\sup\lambda_i(M,g)=+\infty,$
where supremum is taken over the space
of Riemannian metrics $g$ on this surface $M.$

But if a functional is invariant under 
rescaling transformations $g\mapsto tg$
then the question about its supremum becomes interesting.
There are several examples of invariant
functionals such as ratios of eigenvalues etc.
We study an invariant functional
$$
\Lambda_i(M,g)=\lambda_i(M,g)\Area(M,g).
$$

It turns out that the question
about the supremum $\sup\Lambda_i(M,g)$
of the functional $\Lambda_i(M,g)$ over the space
of Riemannian metrics $g$ on a fixed surface
$M$ is very difficult and only few results are known.

It is known that functionals $\Lambda_i(M,g)$ are bounded from above.
It was proven in the paper~\cite{Yang-Yau} by Yang and Yau that
for an orientable surface $M$ of genus $\gamma$ the following
inequality holds,
$$
\Lambda_1(M,g)\leqslant 8\pi(\gamma+1).
$$
Korevaar proved in the paper~\cite{Korevaar}
that there exists a constant $C$ such that for
any $i>0$ and any compact surface $M$ of genus $\gamma$
the functional $\Lambda_i(M,g)$ is bounded,
$$
\Lambda_i(M,g)\leqslant C(\gamma+1)i.
$$
However, Colbois and Dodziuk proved in the
paper~\cite{Colbois-Dodziuk}
that for a manifold $M$ of dimension $\dim M\geqslant 3$
the functional $\Lambda_i(M,g)$ 
is not bounded on the space of Riemannian metrics $g$ on $M.$

The functional $\Lambda_i(M,g)$ depends continuously
on the metric $g,$ but this functional
is not differentiable. However, Berger proved in 
the paper~\cite{Berger} that for analytic
deformations $g_t$ the left and right derivatives
of the functional $\Lambda_i(M,g_t)$ with respect to $t$ exist.
This is motivation for the  following definition,
see the paper~\cite{Nadirashvili1}
by Nadirashvili and the paper~\cite{ElSoufi-Ilias1}
by El Soufi and Ilias.

\begin{Definition}
A Riemannian metric $g$ on  a closed surface
$M$ is called extremal metric for the
functional $\Lambda_i(M,g)$ if for any analytic deformation
$g_t$ such that $g_0=g$ the following inequality holds,
$$
\frac{d}{dt}\Lambda_i(M,g_t)%
\left.\vphantom{\raisebox{-0.5em}{.}}\right|_{t=0+}\leqslant0%
\leqslant\frac{d}{dt}\Lambda_i(M,g_t)%
\left.\vphantom{\raisebox{-0.5em}{.}}\right|_{t=0-}.
$$
\end{Definition}

The list of surfaces $M$ and values of index $i$
such that maximal or at least extremal metrics
for the functional $\Lambda_i(M,g)$ are known is quite short.
Let us denote the Klein bottle by $\mathbb{K}$
and the surface of genus $2$ by $\Sigma_2.$

\begin{enumerate}
\item Complete answers for $\Lambda_1(\mathbb{S}^2,g),$
$\Lambda_2(\mathbb{S}^2,g),$ $\Lambda_1(\mathbb{R}P^2,g),$
$\Lambda_1(\mathbb{T}^2,g)$ and $\Lambda_1(\mathbb{K},g)$
are known, and a complete answer for $\Lambda_1(\Sigma_2)$
is conjectured.
\item There are also some partial results concerning
$\Lambda_i(\mathbb{T}^2,g)$ and $\Lambda_i(\mathbb{K},g)$
for $i>1.$
\end{enumerate}

Let us describe these results in details.
Hersch proved in the paper~\cite{Hersch}
that $\sup\Lambda_1(\mathbb{S}^2,g)=8\pi$ and the maximum
is reached on the canonical metric on $\mathbb{S}^2.$
This metric is the unique extremal metric
for the functional $\Lambda_1(\mathbb{S}^2,g).$
Nadirashvili proved in the paper~\cite{Nadirashvili2}
that $\sup\Lambda_2(\mathbb{S}^2,g)=16\pi$
and maximum is reached on a singular metric which can be obtained
as the metric on the union of two spheres of equal radius with
canonical metric glued together.

Li and Yau proved in the paper~\cite{Li-Yau}
that $\sup\Lambda_1(\mathbb{R}P^2,g)=12\pi$
and the maximum is reached on the canonical 
metric on $\mathbb{R}P^2.$ This metric is the 
unique extremal metric for the functional
$\Lambda_1(\mathbb{R}P^2,g).$

Nadirashvili
proved in the paper~\cite{Nadirashvili1}
that $\sup\Lambda_1(\mathbb{T}^2,g)=\frac{8\pi^2}{\sqrt{3}}$
and the maximum is reached on the flat equilateral torus. El~Soufi
and Ilias proved in the paper~\cite{ElSoufi-Ilias1} 
that the only extremal
metric for $\Lambda_1(\mathbb{T}^2,g)$ different from the maximal one
is the metric on the Clifford torus.

Jakobson, Nadirashvili and 
Polterovich proved in 
the paper~\cite{Jakobson-Nadirashvili-Polterovich}
that the metric on a Klein bottle realized
as the Lawson bipolar surface $\tilde{\tau}_{3,1}$
is extremal. El Soufi, Giacomini and 
Jazar proved in the
paper~\cite{ElSoufi-Giacomini-Jazar}
that this metric is the unique extremal metric and the maximal
one.

It is shown by Jakobson, Levitin, Nadirashvili, Nigam and Polterovich
in the paper~\cite{JLNNP} using a combination of analytic and numerical tools
that the maximal metric for the first eigenvalue on 
the surface of genus two is the metric on the Bolza surface $\mathcal P$
induced from the canonical metric on the sphere
using the standard covering ${\mathcal P}\longrightarrow\mathbb{S}^2.$ 
In fact, the authors state 
this result as a conjecture, because a part of the argument 
is based on a numerical calculation.

Let $m,k\in\mathbb{N},$ $0<k<m,$ $(m,k)=1.$
Lapointe studied in the paper~\cite{Lapointe}
Lawson bipolar surfaces $\tilde{\tau}_{m,k}$ 
and proved the following result. 
\begin{enumerate}
\item If $mk\equiv 0 \mod 2$ then
$\tilde{\tau}_{m,k}$ is a torus and it carries an extremal metric
for $\Lambda_{4m-2}(\mathbb{T}^2,g).$
\item If $mk\equiv 1 \mod 4$ then
$\tilde{\tau}_{m,k}$ is a torus and it carries an extremal metric
for $\Lambda_{2m-2}(\mathbb{T}^2,g).$
\item If $mk\equiv 3 \mod 4$ then
$\tilde{\tau}_{m,k}$ is a Klein bottle and it carries an extremal metric
for $\Lambda_{m-2}(\mathbb{K},g).$
\end{enumerate}

The example of $\Lambda_2(\mathbb{S}^2,g)$ shows 
that a maximal metric could be {\em singular}.
All above mentioned results were obtained by different
complicated techniques invented ad hoc
for different surfaces and eigenvalue numbers.

A general approach for finding {\em smooth} extremal
metrics became possible due to
a result by  El Soufi and Ilias published
in the paper~\cite{ElSoufi-Ilias2}. They proved
that a metric on a compact immersed
minimal surface in $\mathbb{S}^n$
is extremal for some functional $\Lambda_i,$
but we should emphasize that the value of 
the corresponding index $i$ is unknown
(see exact statement
in Proposition~\ref{minimal-extremal}).

Several immersed minimal
surfaces in $\mathbb{S}^n$ are known.
Therefore, at least theoretically
we can start with such a surface and
find the corresponding functional $\Lambda_j$
such that the metric on the surface is extremal
for this functional. However,  it turns out
that finding the value of the index $j$
is {\em a difficult problem.}

The author proposed to find this index $j$
by a careful analysis of eigenvalue
structure. For the first time this approach
was successfully realized
in the  recent paper~\cite{Penskoi1}
where Lawson tau-surfaces $\tau_{m,k}$
were studied. The Lawson tau-surfaces are tori
or Klein bottles minimally immersed in $\mathbb{S}^3.$
Here $m,k\geqslant1$ are integers and we can always 
assume that $(m,k)=1.$
The following result is proved.
\begin{enumerate}
\item If $\tau_{m,k}$ is a Lawson torus, i.e. $m,k\equiv 1 \mod 2,$
then it carries an extremal metric
for the functional $\Lambda_j(\mathbb{T}^2,g),$
where
$j=2\left[\frac{\sqrt{m^2+k^2}}{2}\right]+m+k-1.$
\item If $\tau_{m,k}$ is a Lawson Klein bottle, i.e. 
we can assume that $m\equiv 0 \mod 2,$ $k\equiv 1 \mod 2,$
then it carries an extremal metric
for the functional $\Lambda_j(\mathbb{K},g),$
where 
$j=2\left[\frac{\sqrt{m^2+k^2}}{2}\right]+m+k-1.$
\end{enumerate}
The proof is quite complicated and became possible
only because it turned out that the spectral problem
for the Laplace-Beltrami operator on a Lawson
tau-surface could be reduced to the classical
Lam\'e equation.

The present paper continue the same approach but
in quite a different setting. We are going to 
investigate Otsuki tori. The definition of the Otsuki
tori is quite complicated and we postpone it
till Section~\ref{otsuki}. At this moment it
is sufficient to know that 
for each rational number $\frac{p}{q}$ such that $p,q>0,$
$(p,q)=1$ and $\frac{1}{2}<\frac{p}{q}<\frac{\sqrt{2}}{2}$
there exists an immersed minimal two-dimensional
torus in $\mathbb{S}^3$ called Otsuki torus.
Let us denote it by $O_\frac{p}{q}.$

While the Lawson tau-surfaces
are defined by a simple explicit parameterization,
the Otsuki tori are defined in quite an implicit way.
While in the case of the Lawson tau-surfaces we solve
the problem by reducing to the Lam\'e equation,
in the case of the Otsuki tori we cannot even write
down the Laplace-Beltrami operator explicitly.
However, we prove the following
Theorem.

\begin{Theorem} The metric on an Otsuki
torus $O_\frac{p}{q}\subset\mathbb{S}^3$
is extremal for the functional
$\Lambda_{2p-1}(\mathbb{T}^2,g).$ 
\end{Theorem}

Let us make the following observations.

\begin{enumerate}
\item The Otsuki tori provide an extremal metric for the
eigenvalues $\Lambda_j$ with {\em any} odd $j\geqslant3.$
In comparison, the bipolar Lawson tori $\tilde{\tau}_{m,k}$
provide extremal metrics for eigenvalues $\Lambda_j$ with some
even $j,$ and the Lawson tori $\tau_{m,k}$ provide extremal metrics
for eigenvalues $\Lambda_j$
with {\em some} odd $j,$ namely $\Lambda_1,$ $\Lambda_{5},$
$\Lambda_{9},$ $\Lambda_{11},$ $\Lambda_{13}$ etc.
\item In particular, the Otsuki torus $O_\frac{2}{3}$
provides an extremal metric for 
$\Lambda_3(\mathbb{T}^2,g).$ Except $\Lambda_1,$
it is the eigenvalue $\Lambda_j$ with the smallest
$j$ such that we know an extremal metric. The next
known result is the functional $\Lambda_5(\mathbb{T}^2)$
with a known extremal metric on the Clifford torus $\tau_{3,1}.$
\end{enumerate}

The paper is organized in the following way. We define the
Otsuki tori in Section~\ref{otsuki} and explain
the relation between minimal surfaces in spheres
and extremal metrics in Section~\ref{minimal}.
Then we prove the Theorem and its corollaries 
in Section~\ref{proof-section}.

\section{Otsuki tori}\label{otsuki}

In the paper~\cite{otsuki} Tominosuke Otsuki
introduced a family of minimal hypersurfaces in spheres.
We are interested only in the case of $\mathbb{S}^3.$
In this particular case the introduced hypersurfaces are 
two-dimensional tori and we call them Otsuki tori. The
original definition by Otsuki is complicated and we will
rather explain how to define the Otsuki tori as a part of
a general construction by Hsiang and Lawson of
minimal submanifolds of low cohomogeneity proposed
in the paper~\cite{hsiang-lawson}.

\subsection{Minimal submanifolds of low cohomogeneity.}

Let $M$ be a Riemannian manifold with
a metric $g'$ and $I(M)$ its
full isometry group. Let $G\subset I(M)$
be an isometry group. The space of orbits $M/G$
is naturally a differentiable stratified set.
Let us denote by $\pi$ the natural
projection $\pi:M\longrightarrow M/G.$

The union $M^*$ of all orbits of principal type
is an open dense submanifold of~$M.$
The subset $M^*/G$ of $M/G$ is a manifold
carrying a natural Riemannian structure $g$
defined by the formula $g(X,Y)=g'(X',Y'),$
where $X,Y$ are tangent vectors at a point
$x\in M^*/G$ and $X',Y'$ are tangent vectors at 
a point $x'\in \pi^{-1}(x)\subset M^*$
such that $d_{x'}\pi(X')=X,$ $d_{x'}\pi(Y')=Y$
and $X',$ $Y'$ are orthogonal to the orbit
$\pi^{-1}(x).$

Let us define a volume function $V:M/G\longrightarrow\mathbb{R}$
in the following way. 
If $x\in M^*/G$ then $V(x)=Vol(\pi^{-1}(x)),$
for other points we define $V$ by continuity.
As a result, $V$ is differentiable on $M^*/G$
and continuous on $M/G.$

Let $f:N\longrightarrow M$ be a $G$-invariant
submanifold, i.e. $G$ acts on $N$ and $f$ commutes with 
the actions of $G$ on $N$ and $M.$

\begin{Definition}[see the paper~\cite{hsiang-lawson}] A cohomogeneity
of a $G$-invariant submanifold $f:N\longrightarrow M$ 
in $M$ is the integer $\dim N-\nu,$
where $\dim N$ is the dimension of $N$ and
$\nu$ is the common dimension of the principal orbits.
\end{Definition}

Let us define for each integer $k\geqslant 1$ a metric
\begin{equation}\label{gk}
g_k=V^\frac{2}{k}g.
\end{equation}

The following Proposition (the reduction theorem) is
an important tool for constructing minimal manifolds.

\begin{Proposition}[Theorem 2 from the paper~\cite{hsiang-lawson}]\label{G-reduction}
Let $f:N\longrightarrow M$ be a $G$-invariant submanifold
of cohomogeneity $k,$ and let $M/G$ be given the metric
$g_k$ defined by formula~\eqref{gk}. Then
$f:N\longrightarrow M$ is minimal is and only if
$\bar{f}:N^*/G\longrightarrow M^*/G$
is minimal.
\end{Proposition}

\subsection{Otsuki tori as minimal $\SO(2)$-invariant
submanifolds of cohomogeneity $1$ in $\mathbb{S}^3$}

Let us now apply this theory in the particular case
of $M=\mathbb{S}^3$ and $G=\SO(2).$ Let $x,y,z,t$
be coordinates in $\mathbb{R}^4$ and $\mathbb{S}^3$
be defined by the equation $x^2+y^2+z^2+t^2=1.$
The action of $\SO(2)$ is given by the formula
$$
\alpha\cdot(x,y,z,t)=%
(\cos\alpha x+\sin\alpha y,-\sin\alpha x+\cos\alpha y,z,t),
$$
where $\alpha\in[0,2\pi)$ is a coordinate on $\SO(2).$
The orbits of principal type are circles of radius
$q=\sqrt{x^2+y^2}>0,$ exceptional orbits
are points. The space of orbits 
$\mathbb{S}^3/\SO(2)$ is the closed half-sphere
$\mathbb{S}^2_+$
$$
q^2+z^2+t^2=1,\quad q\geqslant0,
$$
a point $(q,z,t)\in\mathbb{S}^2_+$
corresponds to an orbit
$(q\cos\alpha,q\sin\alpha,z,t)\in\mathbb{S}^3.$
It is easy to see that $(\mathbb{S}^3)^*$
is the open submanifold of $\mathbb{S}^3$
given by the inequality $x^2+y^2>0$
and its image $(\mathbb{S}^3)^*/\SO(2)$
in the space of orbits is the open
half-sphere
$\mathbb{S}^2_{>0}$
$$
q^2+z^2+t^2=1,\quad q>0.
$$

The volume function is given by the formula
$V(q,z,t)=2\pi q$ since it is the length
of a circle of radius $q.$

Let us introduce spherical coordinates 
$0\leqslant\varphi\leqslant\frac{\pi}{2},$
$0\leqslant\theta<2\pi$
in the space of orbits
$\mathbb{S}^3/\SO(2)=\mathbb{S}^2_+,$
$$
\left\{\begin{array}{l}
q=\sin\varphi,\\
z=\cos\varphi\cos\theta,\\
t=\cos\varphi\sin\theta.
\end{array}\right.
$$
In these terms the volume function is given
by the formula $V(\varphi,\theta)=2\pi\sin\varphi.$

Adding the coordinate $\alpha\in[0,2\pi)$ along 
the orbits, we obtain coordinates
on $\mathbb{S}^3$ convenient from the point
of view of the group action,
\begin{equation}\label{s3-coords}
\left\{\begin{array}{l}
x=\cos\alpha\sin\varphi,\\
y=\sin\alpha\sin\varphi,\\
z=\cos\varphi\cos\theta,\\
t=\cos\varphi\sin\theta.
\end{array}\right.
\end{equation}
Then it is easy to find the metric on $\mathbb{S}^3,$
\begin{equation}\label{gprime}
g'=\sin^2\varphi d\alpha^2+d\varphi^2+\cos^2\varphi d\theta^2,
\end{equation}
and the induced metric on the space of orbits
$\mathbb{S}^3/\SO(2)=\mathbb{S}^2_+,$
$$
g=d\varphi^2+\cos^2\varphi d\theta^2.
$$

Let us now look for minimal $\SO(2)$-invariant
submanifolds $f:N\longrightarrow\mathbb{S}^3$
of cohomogeneity $1.$ Since the principal
orbits are circles, by Proposition~\ref{G-reduction}
these submanifolds correspond to one-dimensional
minimal submanifolds $N^*/\SO(2)$ in
the manifold $(\mathbb{S}^3)^*/\SO(2)$
equipped with the metric
\begin{equation}\label{g1}
g_1=V^2g=4\pi^2\sin^2\varphi(d\varphi^2+\cos^2\varphi d\theta^2).
\end{equation}
But one-dimensional minimal submanifolds are just
closed geodesics. Hence, we obtain the following procedure
providing minimal 2-dimensional manifolds in $\mathbb{S}^3,$
\begin{enumerate}
\item find closed geodesics $\gamma$ in
$(\mathbb{S}^3)^*/\SO(2)=\mathbb{S}^2_{>0}$
equipped with the metric~\eqref{g1},
\item find their preimages $N=\pi^{-1}(\gamma)$
in $\mathbb{S}^3.$
\end{enumerate}

We should remark that one has to be careful
with the word "closed" since there exist geodesics
on $(\mathbb{S}^3)^*/\SO(2)=\mathbb{S}^2_{>0}$
with endpoints in the complement
to $(\mathbb{S}^3)^*/\SO(2)$ in $\mathbb{S}^3/\SO(2),$
i.e. in the equator $q=0.$ Such geodesics are also
closed in some sense, but we are not interested in them
and refer the reader to the paper~\cite{hsiang-lawson} for
more details. We consider only geodesics in
$(\mathbb{S}^3)^*/\SO(2)=\mathbb{S}^2_{>0}$
closed in the usual sense.

Since a closed geodesic $\gamma$ is homeomorphic to 
the circle $\mathbb{S}^1$ and the group $\SO(2)$
is also homeomorphic to the circle $\mathbb{S}^1,$
the obtained minimal surface $N=\pi^{-1}(\gamma)$
is an immersed $\SO(2)$-invariant minimal 
two-dimensional torus in $\mathbb{S}^3.$

\begin{Definition}\label{otsuki-definition}
An immersed minimal $\SO(2)$-invariant
two-dimensional torus in $\mathbb{S}^3$
such that its image by the projection
$\pi:\mathbb{S}^3\longrightarrow\mathbb{S}^3/\SO(2)$
is a closed geodesics in
$(\mathbb{S}^3)^*/\SO(2)$
equipped with the metric~\eqref{g1}
is called an Otsuki torus. 
\end{Definition}

Obviously, this definition is not constructive.
Unfortunately, it is not possible to define
Otsuki tori by an equation or by an explicit
parameterization. However, Otsuki tori could
be described in some terms.

\begin{Proposition}[see the papers~\cite{hsiang-lawson,otsuki}]\label{otsuki-pq}
Except one particular case given
by the equation $\varphi=\frac{\pi}{4}$
(this is a Clifford torus),
Otsuki tori are in one-to-one correspondence
with rational numbers $\frac{p}{q}$
such that
$$
\frac{1}{2}<\frac{p}{q}<\frac{\sqrt{2}}{2},\quad p,q>0,\quad (p,q)=1.
$$ 
\end{Proposition}

\begin{Definition} By $O_\frac{p}{q}$ we denote
the Otsuki torus corresponding to $\frac{p}{q}.$
\end{Definition}

Then it is convenient to amend Definition~\ref{otsuki-definition}
and reserve the term "Otsuki tori" only for tori $O_\frac{p}{q}.$

\noindent{\bf Proof} of Proposition~\ref{otsuki-pq} 
could be found, in fact, in the paper~\cite{otsuki}.
We write "in fact" because the statement and the proof
are divided there in several pieces dispersed through
a lengthy paper. We expose here the proof
in more organized way.

In order to prove Proposition~\ref{otsuki-pq}, one should investigate
closed geodesics on $(\mathbb{S}^3)^*/\SO(2)$
equipped with the metric~\eqref{g1}.
Let us use standard notation
$$
E=4\pi^2\sin^2\!\varphi,\quad G=4\pi^2\sin^2\!\varphi\cos^2\!\varphi
$$
for the coefficients of the metric $g_1.$

The equation of geodesics for $\ddot\theta$ reads
$$
\ddot\theta+\frac{1}{G}\frac{\partial G}{\partial\varphi}\dot\varphi\dot\theta=0%
\Longleftrightarrow \frac{d}{dt}(G\dot\theta)=0.
$$
Hence, $c=G\dot\theta$ is an integral of motion and
\begin{equation}\label{theta-dot}
\dot\theta=\frac{c}{G}.
\end{equation}

As we know, the velocity vector of a geodesic is of
constant length. Let us put it equal to $1.$ Hence, we have
\begin{equation}\label{phi-dot}
E\dot\varphi^2+G\dot\theta^2=1\Longleftrightarrow%
\dot\varphi^2=\frac{G-c^2}{EG}.
\end{equation}
This implies $G-c^2\geqslant0$ and
$G=c^2$ if and only if $\dot\varphi=0.$
Let us remark that 
on the interval $0\leqslant\varphi\leqslant\frac{\pi}{2}$
the function $G=4\pi^2\sin^2\!\varphi\cos^2\!\varphi$ 
is strictly increasing on $0<\varphi<\frac{\pi}{4},$
strictly decreasing on $\frac{\pi}{4}<\varphi<\frac{\pi}{2},$
and is invariant with respect to the transformation 
$\varphi\mapsto\frac{\pi}{2}-\varphi.$ 

Let us remark that the point $\varphi=0$
does not belong to $(\mathbb{S}^3)^*/\SO(2).$
Hence, a closed geodesic does not pass through this point
and there exist a minimal value $a$ of the coordinate
$\varphi$ on this geodesic. It follows that
$c=\sqrt{G}|_{\varphi=a}=\pm2\pi\sin a\cos a$
and the geodesic is situated in the annulus 
$a\leqslant\varphi\leqslant\frac{\pi}{2}-a.$
Let us choose the parameter $t$ in such a way that
$\varphi(0)=a,$ hence $\dot\varphi(0)=0.$
We can suppose without loss of generality 
that $\theta(0)=0,$ $\dot\theta(0)>0.$

We should remark that there exist a
particular case $a=\frac{\pi}{4}.$ In this
case $\frac{\pi}{2}-a=a.$ It follows that 
$\varphi(t)\equiv a=\frac{\pi}{4}.$ This is
the exceptional case mentioned in the statement
of Proposition~\ref{otsuki-pq}. One obtain immediately
that this is a Clifford torus in $\mathbb{S}^3$
with a well-known parameterization
$$
\left\{\begin{array}{r}
x=\frac{1}{\sqrt{2}}\cos\alpha,\\
y=\frac{1}{\sqrt{2}}\sin\alpha,\\
z=\frac{1}{\sqrt{2}}\cos\theta,\\
t=\frac{1}{\sqrt{2}}\sin\theta.
\end{array}
\right.
$$
We will ignore this case and suppose that $a\ne\frac{\pi}{4}.$

Equations~\eqref{theta-dot}, \eqref{phi-dot} imply
that
$$
\frac{d\varphi}{d\theta}=\pm\frac{\sqrt{G}\sqrt{G-c^2}}{c\sqrt{E}}=%
\pm\cos\varphi\sqrt{\frac{\sin^2\!\varphi\cos^2\!\varphi}{\sin^2\!a\cos^2\!a}-1}.
$$
The right hand side of this equation becomes $0$ only
at $\varphi=a$ and $\varphi=\frac{\pi}{2}-a$ and it is
clear that $\varphi$ is increasing from $a$ to $\frac{\pi}{2}-a,$
then decreasing from $\frac{\pi}{2}-a$ to $a,$ and then this
pattern repeats.

The crucial question is what is the difference $\Omega(a)$
between the value of $\theta,$ corresponding to $\varphi=a,$
and the nearest value of $\theta,$ 
corresponding to $\varphi=\frac{\pi}{2}-a.$
We have
$$
\frac{d\theta}{d\varphi}=%
\pm\frac{1}{\cos\varphi\sqrt{%
\frac{\sin^2\!\varphi\cos^2\!\varphi}{\sin^2\!a\cos^2\!a}-1}}
$$
and
$$
\Omega(a)=\int_{a}^{\frac{\pi}{2}-a}%
\frac{d\varphi}{\cos\varphi\sqrt{%
\frac{\sin^2\!\varphi\cos^2\!\varphi}{\sin^2\!a\cos^2\!a}-1}}.
$$

It is clear that {\em the geodesic is closed if and only if
$\Omega(a)=\frac{p}{q}\pi$ for some rational number $\frac{p}{q}.$}
Hence, the closed geodesics are in one-to-one correspondence with
rational numbers $\frac{p}{q}$ such that there exist
$a\in(0,\frac{\pi}{4})$ such that $\Omega(a)=\frac{p}{q}\pi.$

The function $\Omega(a)$ is a complicated transcendental function
that could be expressed in terms of the complete elliptic integrals
of first and third kind, see the paper~\cite{hsiang-lawson}. However, one can
prove (see the paper~\cite{otsuki} for a detailed proof) that $\Omega(a)$
possesses the following properties,
\begin{enumerate}
\item $\Omega(a)$ is continuous on $(0,\frac{\pi}{4}]$ and strictly
increasing,
\item $\lim\limits_{a\rightarrow 0+}\Omega(a)=\frac{\pi}{2},$
\item $\Omega(\frac{\pi}{4})=\frac{\sqrt{2}}{2}\pi.$
\end{enumerate}
It follows immediately that the Otsuki tori are in one-to-one
correspondence with rational numbers $\frac{p}{q}$
such that $p,q>0,$ $(p,q)=1$ and
$\frac{\pi}{2}<\frac{p}{q}\pi<\frac{\sqrt{2}}{2}\pi.$
This implies Proposition~\ref{otsuki-pq}. $\Box$

\begin{figure}[h]
\begin{center}
\includegraphics{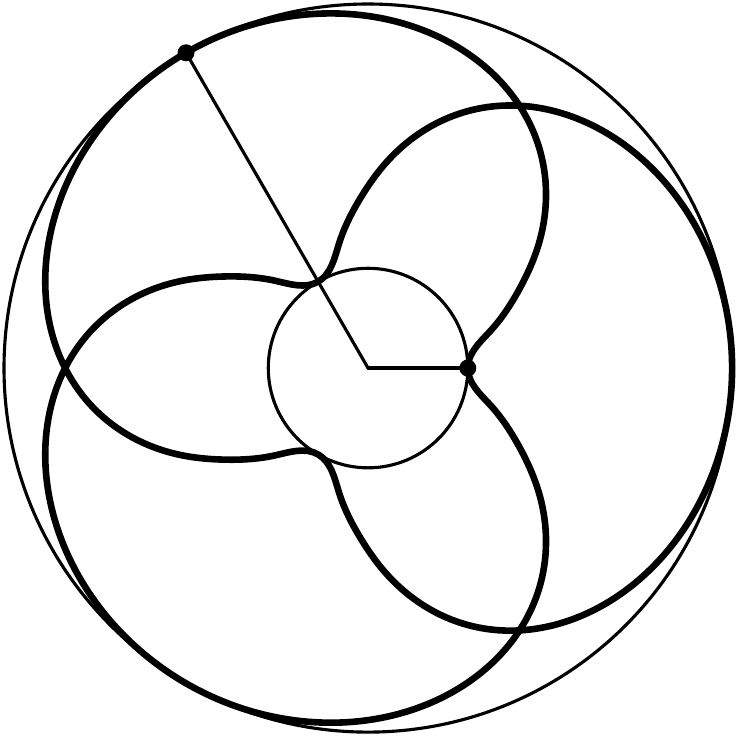}
\end{center}
\caption{The geodesic $\pi(O_\frac{2}{3})$ in
the space of orbits. We use $\varphi$ and $\theta$ as
polar coordinates in the plane.
The circles $\varphi=a$ and $\varphi=\frac{\pi}{2}-a$ are
also drawn. The angle between
the directions pointing towards consecutive minimum and
maximum of $\varphi$ is equal to $\Omega(a)=\frac{2}{3}\pi.$}
\label{o23}
\end{figure}

Let us consider the Otsuki torus $O_\frac{2}{3}$ as
an example. The corresponding value of $a$ such 
that $\Omega(a)=\frac{2}{3}\pi$ could be found
only numerically, $a=0.33787\dots$ The corresponding
closed geodesic $\pi(O_\frac{2}{3})$ is presented
on Fig.~\ref{o23}.

\begin{Proposition}\label{theta-zero}
The functions $\cos\theta$ and $\sin\theta$ have
exactly $2p$ zeroes on a geodesic $\pi(O_\frac{p}{q}).$
\end{Proposition}

\noindent{\bf Proof.} $\Omega(a)=\frac{p}{q}\pi$ is the difference
between values of $\theta$ in consecutive minimum and
maximum of $\varphi,$ see Fig.~\ref{o23} for an example.
It follows that the difference
between values of $\theta$ in two consecutive minimums
is $2\Omega(a)=\frac{p}{q}2\pi.$ Since the geodesic
is closed, there exist smallest possible integers $m>0$ and $n>0$
such that $\frac{p}{q}2\pi n=m2\pi,$ i.e. 
$\frac{p}{q}n=m.$ Since $(p,q)=1,$ we obtain $n=q$ and $m=p.$

This implies that the natural projection
$\pi(O_\frac{p}{q})\stackrel{\theta}\longrightarrow\mathbb{R}/2\pi\mathbb{R}$
is a $p$-fold covering. Since $\cos\theta$ and $\sin\theta$ have
exactly $2$ zeroes on $[0,2\pi),$ we obtain the Proposition. $\Box$

\section{Minimal submanifolds of a sphere
and extremal spectral property of their metrics}\label{minimal}

Let us recall two important results about
minimal submanifolds of a sphere.
Let $N$ be a $d$-dimensional minimal
submanifold of the sphere $\mathbb{S}^n\subset\mathbb{R}^{n+1}$
of radius~$R.$ Let $\Delta$ be the Laplace-Beltrami operator
on $N$ equipped with the induced metric.

The first result is a classical one. 
Its proof can be found e.g. in the book~\cite{Kobayashi-Nomizu}.

\begin{Proposition}\label{minimal-eigenfunction}
The restrictions $x^1|_N,\dots,x^{n+1}|_N$ on $N$
of the standard coordinate functions of $\mathbb{R}^{n+1}$ 
are eigenfunctions of $\Delta$ with eigenvalue $\frac{d}{R^2}.$
\end{Proposition}

Let us numerate the eigenvalues of $\Delta$ as in formula~\eqref{eigenvalues}
counting them with multiplicities
$$
0=\lambda_0<\lambda_1\leqslant\lambda_2\leqslant\dots\leqslant\lambda_i\leqslant\dots
$$
Proposition~\ref{minimal-eigenfunction} implies that there exists
at least one index $i$ such that $\lambda_i=\frac{d}{R^2}.$
Let $j$ denotes the minimal number $i$ such that $\lambda_i=\frac{d}{R^2}.$
Let us introduce the eigenvalues counting function
$$
N(\lambda)=\#\{\lambda_i|\lambda_i<\lambda\}.
$$
We see that $j=N(\frac{d}{R^2}).$ Remark that we count the
eigenvalues starting from $\lambda_0=0.$

The second result was published in 2008 by 
El~Soufi and Ilias in the paper~\cite{ElSoufi-Ilias2}.
It could be written in the following form.

\begin{Proposition}\label{minimal-extremal}
The metric induced on $N$ by immersion
$N\subset\mathbb{S}^n$ is an extremal metric
for the functional $\Lambda_{N\left(\frac{d}{R^2}\right)}(N,g).$
\end{Proposition}

We should remark that isometric immersions by eigenfunctions
of the Laplace-Beltrami operator were studied
in the paper~\cite{Takahashi} by Takahashi. The results
of Takahashi are used by El~Soufi and Ilias in the 
paper~\cite{ElSoufi-Ilias2}.

Proposition~\ref{minimal-extremal} implies immediately the following one.

\begin{Proposition}\label{reduction}
The  metric induced on an Otsuki torus $O_\frac{p}{q}$
by its immersion $O_\frac{p}{q}\subset\mathbb{S}^3$
is an extremal metric for the functional 
$\Lambda_{N(2)}(\mathbb{T}^2,g).$
\end{Proposition}

\noindent{\bf Proof.} Since an Otsuki torus 
is a complete minimal surface in the sphere $\mathbb{S}^3$
of radius $1,$
the statement follows immediately from  
Proposition~\ref{minimal-extremal}, where $R=1$ and $d=2.$
$\Box$

Proposition~\ref{reduction} is crucial for this paper. It reduces
investigation of extremal spectral properties of the Otsuki tori
to counting eigenvalues $\lambda_i$ of the Laplace-Beltrami
operator such that $\lambda_i<2.$

\section{Proof of the Theorem and its corollaries}\label{proof-section}

\subsection{Eigenvalues of the Laplace-Beltrami
operator on the Otsuki tori and auxiliary periodic
Sturm-Liouville problem}

Let us now reduce the counting eigenvalues problem
for the Laplace-Beltrami operator on the Otsuki tori
to the same problem for the periodic Sturm-Liouville
problem following the same procedure as in the  
paper~\cite{Penskoi1}.

Let us consider an Otsuki torus~$O_\frac{p}{q}.$
As coordinates we take the parameter 
$\alpha\in[0,2\pi)$ on $\SO(2)$ and a natural
parameter $t$ on the geodesic $\pi(O_\frac{p}{q})$
such that $t$ takes values in the interval $[0,t_0),$
where $t_0$ is the length of the closed geodesic
$\pi(O_\frac{p}{q})$ with respect to
the metric $g_1$ given by formula~\eqref{g1}.
Let $\varphi(t)$ and $\theta(t)$ be the parameterization
of the geodesic $\pi(O_\frac{p}{q})$ by $t.$ We should
emphasize that it is impossible to find 
$\varphi(t)$ and $\theta(t)$ explicitly,
but {\em it turns out that it does not
matter.}

\begin{Proposition}
Let $O_\frac{p}{q}$ be an Otsuki
torus parameterized by $\alpha$ and $t$
as described above.
Then the metric induced by the immersion
$O_\frac{p}{q}\subset\mathbb{S}^3$
is equal to
\begin{equation}\label{o-metric}
\sin^2\!\varphi(t)d\alpha^2+\frac{1}{4\pi^2\sin^2\!\varphi(t)}dt^2
\end{equation}
and the Laplace-Beltrami operator is given by 
the following formula,
\begin{equation}\label{Delta}
\Delta f=-\frac{1}{\sin^2\!\varphi(t)}\frac{\partial^2 f}{\partial\alpha^2}-%
\frac{\partial}{\partial t}%
\left(4\pi^2\sin^2\!\varphi(t)\frac{\partial f}{\partial t}\right).
\end{equation}
\end{Proposition}

\noindent{\bf Proof.} The metric $g'$ on the
sphere $\mathbb{S}^3$ is given by formula~\eqref{gprime}.
Hence, the metric on an Otsuki torus $O_\frac{p}{q}$
is given by the formula 
$$
\sin^2\!\varphi(t)d\alpha^2+(\dot\varphi^2(t)+\cos^2\!\varphi(t)%
\dot\theta^2(t))dt^2.
$$
But it follows from equation~\eqref{phi-dot} that
$$
\dot\varphi^2(t)+\cos^2\!\varphi(t)\dot\theta^2(t)=\frac{1}{4\pi^2\sin^2\!\varphi(t)}
$$
and this implies formula~\eqref{o-metric}.
Formula~\eqref{Delta} is obtained
by a direct calculation. $\Box$

Counting eigenvalues is known to be a difficult problem.
Fortunately, we can reduce this problem to a one-dimensional
one. Let us consider a family of periodic Sturm-Liouville
problems indexed by a parameter $l,$
\begin{gather}
\frac{d}{dt}\left(4\pi^2\sin^2\!\varphi(t)\frac{dh(t)}{d t}\right)+%
\left(\lambda-\frac{l^2}{\sin^2\!\varphi(t)}\right)h(t)=0,\label{Sturm-Liouville}\\
h(t+t_0)\equiv h(t).\label{periodic-h}
\end{gather}
Let $h(l,t)$ denote a solution of
the periodic Sturm-Liouville 
problem~\eqref{Sturm-Liouville}, \eqref{periodic-h}.
We consider $l$ as a parameter in equation~\eqref{Sturm-Liouville}
and in its solution $h(l,t).$ For example, we assume
$l$ fixed and consider $t$ as an independent 
variable when we discuss
zeroes of the function $h(l,t).$

Equation~\eqref{Sturm-Liouville}
is written in the standard form of a Sturm-Liouville problem,
and the following classical result holds, 
see e.g. the book~\cite{Coddington-Levinson}.

\begin{Proposition}\label{SL-properties}
There are an infinite number of eigenvalues $\lambda_i(l)$ 
of the periodic Sturm-Liouville problem~\eqref{Sturm-Liouville},
\eqref{periodic-h}.
Given a fixed integer $l,$
eigenvalues $\lambda_i(l)$  form a sequence such that
$$
\lambda_0(l)<\lambda_1(l)\leqslant\lambda_2(l)<%
\lambda_3(l)\leqslant\lambda_4(l)<%
\lambda_5(l)\leqslant\lambda_6(l)<\dots
$$
For $\lambda=\lambda_0(l)$ there exists a unique
(up to
a multiplication by a non-zero constant)
eigenfunction $h_0(l,t).$ If 
$\lambda_{2i+1}(l)<\lambda_{2i+2}(l)$
for some $i\geqslant0$ then there is a unique (up to
a multiplication by a non-zero constant) eigenfunction 
$h_{2i+1}(l,t)$ with eigenvalue $\lambda=\lambda_{2i+1}(l)$
of multiplicity one
and there is a unique (up to
a multiplication by a non-zero constant) eigenfunction 
$h_{2i+2}(l,t)$ with eigenvalue $\lambda=\lambda_{2i+2}(l)$
of multiplicity one.
If  $\lambda_{2i+1}(l)=\lambda_{2i+2}(l)$
then there are two independent eigenfunctions
$h_{2i+1}(l,t)$ and $h_{2i+2}(l,t)$ 
with eigenvalue $\lambda=\lambda_{2i+1}(l)=\lambda_{2i+1}(l)$
of multiplicity two.

The eigenfunction $h_0(l,t)$ has no zeroes on $[0,t_0).$
The eigenfunctions 
$h_{2i+1}(l,t)$ and $h_{2i+2}(l,t),$
$i\geqslant0,$ each  
have exactly $2i+2$ zeroes on $[0,t_0).$
\end{Proposition}

Let us now explain a relation between the family
of periodic Sturm-Liouville problems~\eqref{Sturm-Liouville},
\eqref{periodic-h} and the spectral problem for the
Laplace-Beltrami operator $\Delta$ on an Otsuki torus $O_\frac{p}{q}.$

\begin{Proposition}\label{reduction-proposition}
A number $\lambda$ is an eigenvalue of $\Delta$
if and only if there exists a non-negative
integer $l$ and an eigenvalue $\lambda_i(l)$
of the periodic Sturm-Liouville problem~\eqref{Sturm-Liouville},
\eqref{periodic-h} with parameter $l$ such that $\lambda_i(l)=\lambda.$

The eigenspace of the Laplace-Beltrami
operator $\Delta$ with eigenvalue $\lambda$
has a basis consisting of functions of the form
\begin{equation}\label{h-cos}
h_i(l,t)\cos(l\alpha),
\end{equation}
where $l=0,1,2,\dots$ and there exists $i$ such that $\lambda_i(l)=\lambda,$
and
\begin{equation}\label{h-sin}
h_i(l,t)\sin(l\alpha),
\end{equation}
where $l=1,2,\dots$ and there exists $i$ such that $\lambda_i(l)=\lambda.$
\end{Proposition}

\noindent{\bf Proof.}
Let us remark that $\Delta$ commutes 
with $\frac{\partial}{\partial\alpha}.$ It follows
that $\Delta$ has a basis
of eigenfunctions of the form $h(l,t)\cos(l\alpha)$
and $h(l,t)\sin(l\alpha)$. Substituting 
these eigenfunctions
into the formula $\Delta f=\lambda f,$ we obtain
equation~\eqref{Sturm-Liouville}. 
However, these solutions should be invariant under
transformations
$$
(\alpha,t)\mapsto(\alpha+2\pi,t),\quad(\alpha,t)\mapsto(\alpha,t+t_0).
$$
This condition implies the condition $l\in\mathbb{Z}$
and the periodicity condition~\eqref{periodic-h}.
It remains only to apply
Proposition~\ref{SL-properties} in order
to describe the basis.
$\Box$

It is easy now to establish a relation between
the multiplicities of the eigenvalues of the 
Laplace-Beltrami operator $\Delta$ on an Otsuki torus
and the eigenvalues $\lambda_i(l)$
of the periodic Sturm-Liouville 
problem~\eqref{Sturm-Liouville},~\eqref{periodic-h}.
This relation permits us to express the quantity $N(2)$
in terms of the eigenvalues $\lambda_i(l).$

\begin{Proposition}\label{how-to-count}
Let $O_\frac{p}{q}$ be an Otsuki torus.
Then 
\begin{equation}\label{N2torus}
N(2)=\#\{\lambda_i(0)|\lambda_i(0)<2\}+%
2\#\{\lambda_i(l)|\lambda_i(l)<2,l>0,l\in\mathbb{Z}\}.
\end{equation}
\end{Proposition}

\noindent{\bf Proof.} 
The eigenvalue $\lambda_i(0)$ gives
exactly one basis eigenfunction 
$$
h_i(0,t)\cos(0\alpha)=h_i(0,t)
$$
of the operator~$\Delta.$
It follows that
each $\lambda_i(0)$ corresponds to one eigenvalue
$\lambda_q=\lambda_i(0)$ of the Laplace-Beltrami operator $\Delta.$
The eigenvalue $\lambda_i(l),$ $l>0,$ gives
exactly two basis eigenfunctions 
$$h_i(l,t)\cos(l\alpha)\quad
\mbox{and}\quad h_i(l,t)\sin(l\alpha)
$$ of the Laplace-Beltrami operator~$\Delta.$
It follows that 
for $l>0$ each $\lambda_i(l)$ corresponds to two eigenvalues
$\lambda_q=\lambda_{q+1}=\lambda_i(l)$ of the Laplace-Beltrami operator $\Delta.$
This implies formula~\eqref{N2torus}.
$\Box$

Let us now investigate properties of eigenvalues 
$\lambda_i(l)$ as functions of $l.$
The following Proposition 
(see e.g. Proposition~15 in the paper~\cite{Penskoi1})
is very important for the proof of the Theorem.

\begin{Proposition}\label{lambda-l}
Let us fix $i.$ If $\lambda_i(l)$ has multiplicity $1$
for all $l\in(0,l_1)$ then
$\lambda_i(l)$ is a strictly increasing function
on $(0,l_1)$.
\end{Proposition}

In particular, the following Proposition holds.

\begin{Proposition}\label{lambda0-increasing}
The zeroth eigenvalue $\lambda_0(l)$ is a strictly
increasing function of $l.$
\end{Proposition}

\noindent{\bf Proof} follows immediately from 
Proposition~\ref{lambda-l} because $\lambda_0(l)$
is always of multiplicity one according to
Proposition~\ref{SL-properties}.
$\Box$

\subsection{Proof of the Theorem}

It follows immediately from formula~\eqref{s3-coords}
that the coordinate functions $x,$ $y,$ $z$ and $t$
being restricted on the Otsuki torus $O_\frac{p}{q}$
are given by formulae
$$
\left\{\begin{array}{l}
x|_{O_\frac{p}{q}}=\cos\alpha\sin\varphi(t),\\
y|_{O_\frac{p}{q}}=\sin\alpha\sin\varphi(t),\\
z|_{O_\frac{p}{q}}=\cos\varphi(t)\cos\theta(t),\\
t|_{O_\frac{p}{q}}=\cos\varphi(t)\sin\theta(t).
\end{array}\right.
$$
By Proposition~\ref{minimal-eigenfunction}
these restrictions are eigenfunctions of eigenvalue $2$
of the Laplace-Beltrami operator $\Delta$
on the Otsuki torus $O_\frac{p}{q}.$

The function
$x|_{O_\frac{p}{q}}$ is of
the form~\eqref{h-cos} with $l=1,$
the function $y|_{O_\frac{p}{q}}$
is of the form~\eqref{h-sin}
with $l=1,$ the functions 
$z|_{O_\frac{p}{q}}$ and
$t|_{O_\frac{p}{q}}$ are of the 
form~\eqref{h-cos} with $l=0.$

It follows that $\sin\varphi(t)$
should be an eigenfunction of the
periodic Sturm-Liouville problem~\eqref{Sturm-Liouville},
\eqref{periodic-h} with $l=1$ and eigenvalue $\lambda=2.$
In the same way, it follows that
$\cos\varphi(t)\cos\theta(t)$
and $\cos\varphi(t)\sin\theta(t)$
should be eigenfunctions of the
periodic Sturm-Liouville problem~\eqref{Sturm-Liouville},
\eqref{periodic-h} with $l=0$ and eigenvalue $\lambda=2.$
We should remark that this could be checked
by a direct and very long straightforward calculation.

Let us now remember that on an Otsuki torus $O_\frac{p}{q}$
we have an inequality 
$$
0<a\leqslant\varphi\leqslant\frac{\pi}{2}-a<\frac{\pi}{2}.
$$
This implies that $\sin\varphi(t)$ {\em has no zeroes}.
Hence, by Proposition~\ref{SL-properties},
we have 
$$
h_0(1,t)=\sin\varphi(t)\quad\mbox{and}\quad\lambda_0(1)=2.
$$

By Proposition~\ref{theta-zero} the functions
$\cos\varphi(t)\cos\theta(t)$
and $\cos\varphi(t)\sin\theta(t)$ have exactly $2p$
zeroes. Hence, by Proposition~\ref{SL-properties}
we have 
$$
h_{2p-1}(1,t)=\cos\varphi(t)\cos\theta(t),\quad
h_{2p}(1,t)=\cos\varphi(t)\sin\theta(t)$$
and
$$
\lambda_{2p-1}(0)=\lambda_{2p}(0)=2.
$$

Let us now prove that $\lambda_0(0),\dots,\lambda_{2p-2}(0)$
are the only eigenvalues $\lambda_i(l)$ such that
$\lambda_i(l)<2.$

By Proposition~\ref{SL-properties} we have
$$
\lambda_0(0)<\lambda_1(0)\leqslant\lambda_2(0)%
<\lambda_3(0)\leqslant\lambda_4(0)<\dots%
<\lambda_{2p-3}(0)\leqslant\lambda_{2p-2}(0)%
<\lambda_{2p-1}(0)=2.
$$
It follows that $\lambda_0(0),\dots,\lambda_{2p-2}(0)<2.$

Using Proposition~\ref{SL-properties} again, we have
$$
2=\lambda_{2p-1}(0)=\lambda_{2p}(0)%
<\lambda_{2p+1}(0)\leqslant\lambda_{2p+2}(0)<\dots.
$$
It follows that $\lambda_i(0)\geqslant2$ if $i>2p-2.$

We know that $\lambda_0(1)=2.$ By Proposition~\ref{SL-properties},
$$
2=\lambda_0(1)%
<\lambda_1(1)\leqslant\lambda_2(1)%
<\lambda_3(1)\leqslant\lambda_4(1)<\dots.
$$
It follows that $\lambda_i(1)\geqslant2$ for any $i.$

By Proposition~\ref{lambda0-increasing}, $\lambda_0(l)>\lambda_0(1)=2$
for all $l>1.$ Then by Proposition~~\ref{SL-properties}
for all $l>1$ we have 
$$
2<\lambda_0(l)<\dots<\lambda_{2i+1}(l)\leqslant\lambda_{2i+2}(l)<\dots
$$
Hence, $\lambda_i(l)>2$ for any $i$ and any $l>1.$

By Proposition~\ref{how-to-count},
$$
N(2)=\#\{\lambda_i(0)|\lambda_i(0)<2\}+%
2\#\{\lambda_i(l)|\lambda_i(l)<2,l>0,l\in\mathbb{Z}\}=%
$$
$$
=\#\{\lambda_0(0),\dots,\lambda_{2p-2}(0)\}=2p-1.
$$

The statement of the Theorem follows from
Proposition~\ref{reduction}.
$\Box$

\subsection{Corollaries and remarks}

\begin{Proposition}
The value $\Lambda_{2p-1}(O_\frac{p}{q})$ of the functional
$\Lambda_{2p-1}(\mathbb{T}^2,g)$ on its extremal
metric on an Otsuki torus $O_\frac{p}{q}\subset\mathbb{S}^3$
equals twice the length of geodesic $\pi(O_\frac{p}{q})$
with respect to the metric $g_1$ given by
formula~\eqref{g1}.
\end{Proposition}

\noindent{\bf Proof.} We already know
that the metric on an Otsuki
torus $O_\frac{p}{q}\subset\mathbb{S}^3$
is extremal for the functional
$\Lambda_{2p-1}(\mathbb{T}^2,g).$
It follows immediately from formula~\eqref{o-metric}
that the volume form is equal to $\frac{1}{2\pi}d\alpha dt.$
Then
$$
\Lambda_{2p-1}(O_\frac{p}{q})=%
\lambda_{2p-1}\Area(O_\frac{p}{q})=%
2\int\limits_0^{2\pi}\int\limits_0^{t_0}\frac{1}{2\pi}d\alpha dt=2t_0,
$$
and $t_0$ was defined as
the length of the geodesic $\pi(O_\frac{p}{q})$
with respect to the metric $g_1$ given by
formula~\eqref{g1}.
$\Box$

Unfortunately, this length $t_0$ could be found only
numerically. The result of numerical calculations
for several Otsuki tori is given in Table~\ref{table}.

\begin{table}[h]
\caption{Results of numerical calculations for
several Otsuki tori $O_\frac{p}{q}.$}\label{table}
\centerline{\begin{tabular}{|c|c|c|c|}
\hline
$\vphantom{\Bigl(}\frac{p}{q}$&$a$&$2p-1$&$\Lambda_{2p-1}(O_\frac{p}{q})$\\
\hline
$\vphantom{\Bigl(}\frac{2}{3}$ & $0.3379\dots$ & $3$ & $79.91\dots$ \\
\hline
$\vphantom{\Bigl(}\frac{3}{5}$ & $0.1273\dots$ & $5$ & $127.7\dots$ \\
\hline
$\vphantom{\Bigl(}\frac{4}{7}$ & $0.07526\dots$ & $7$ & $177.2\dots$ \\
\hline
$\vphantom{\Bigl(}\frac{5}{8}$ & $0.1874\dots$ & $9$ & $206.7\dots$ \\
\hline
$\vphantom{\Bigl(}\frac{5}{9}$ & $0.05220\dots$ & $9$ & $227.1\dots$ \\
\hline
\end{tabular}}

\end{table}

As we see, the Otsuki torus $O_\frac{2}{3}$
provides an extremal metric for 
$\Lambda_3(\mathbb{T}^2,g).$ Except $\Lambda_1,$
it is the eigenvalue $\Lambda_j$ with the smallest
$j$ such that we know an extremal metric. The next
known result is the functional $\Lambda_5(\mathbb{T}^2)$
with extremal metric on the Lawson torus $\tau_{3,1},$
see the paper~\cite{Penskoi1}.

\begin{Proposition}
Otsuki tori provide an extremal metric for the
eigenvalues $\Lambda_j$ with {\em any} odd $j\geqslant3.$
\end{Proposition}

\noindent{\bf Proof} follows immediately from
existence of a rational number $\frac{p}{q}$
such that $p,q>0,$ $(p,q)=1,$ $2p-1=j\geqslant3$
and
$$\frac{1}{2}<\frac{p}{q}<\frac{\sqrt{2}}{2}.
$$
$\Box$

In comparison, bipolar Lawson tori $\tilde{\tau}_{m,k}$
provide extremal metrics for eigenvalues $\Lambda_j$ with some
even $j,$ and Lawson tori $\tau_{m,k}$ provide extremal metrics
for eigenvalues $\Lambda_j$
with {\em some} odd $j,$ namely $\Lambda_1,$ $\Lambda_{5},$
$\Lambda_{9},$ $\Lambda_{11},$ $\Lambda_{13}$ etc.

One can ask a natural question: if an Otsuki torus and a Lawson torus
carry metrics extremal for the same eigenvalue, which metric provides
larger $\Lambda$? We do not know the answer since 
$\Lambda_{2p-1}(O_\frac{p}{q})$ could be found only numerically,
at least at this moment. However, in all examples studied numerically
Otsuki tori provided larger $\Lambda.$

For example, we know (see the paper~\cite{Penskoi1}) that the Lawson torus $\tau_{3,1}$ carries 
an extremal metric for $\Lambda_{5}$ and 
$$
\Lambda_{5}(\tau_{3,1})=83.98\dots
$$
The Otsuki torus 
$O_\frac{3}{5}$ 
also carries an extremal metric for $\Lambda_{5}$ and
$$
\Lambda_{5}(O_\frac{3}{5})=127.7\dots
$$

\section*{Acknowledgments}

The author thanks I.~Polterovich and P.~Winternitz for fruitful
discussions.

This work was partially supported
by Russian Federation Government grant no.~2010-220-01-077, 
ag. no.~11.G34.31.0005,
by the Russian Foundation
for Basic Research (grant no.~08-01-00541
and grant no.~11-01-12067-ofi-m-2011),
by the Russian State Programme for the Support of
Leading Scientific Schools (grant no.~5413.2010.1)
and by the Simons-IUM fellowship.

\end{document}